%% file: main.tex
\documentclass{article}
\usepackage{ijcai20}
\pdfoutput=1

\usepackage{times}
\usepackage{soul}
\usepackage{url}
\usepackage[utf8]{inputenc}
\usepackage[small]{caption}
\usepackage{subcaption}
\usepackage{graphicx}
\usepackage{amsmath}
\usepackage{amsthm}
\usepackage{amsfonts}
\usepackage{booktabs}
\usepackage{algorithm}
\usepackage{algorithmic}
\usepackage{wrapfig}
\usepackage{tikz}
\usetikzlibrary{shapes.geometric}
\usetikzlibrary{shapes.multipart}
\usepackage{siunitx}
\usepackage[normalem]{ulem}
\usepackage{mathtools}
\DeclarePairedDelimiter{\ceil}{\lceil}{\rceil}

\newcommand{\oname}{{\sc M-RTRS}}
\newcommand{\name}{{\sc A-RTRS}}
\newcommand{\oracle}{{\sc OA-RTRS}}

\title{Real-Time Dispatching of Large-Scale Ride-Sharing Systems: \\ Integrating Optimization, Machine Learning, and Model Predictive Control}

\author{
Connor Riley
\and
Pascal Van Hentenryck
\and
Enpeng Yuan
\affiliations
Georgia Institute of Technology, Atlanta, USA
\emails
\{ctriley,pvh,eyuan8\}@isye.gatech.edu}

\newcommand{\citep}[1]{\citeauthor{#1} \shortcite{#1}}
\begin{document}

\maketitle

\input{abstract.tex}
\input{introduction.tex}
\input{related-work.tex}
\input{rtdar.tex}
\input{prediction.tex}

\input{vehicle-relocation.tex}

\input{results.tex}

\input{conclusion.tex}

\bibliographystyle{named}
\bibliography{references}

\end{document}

%% file: abstract.tex
\begin{abstract}
  This paper considers the dispatching of large-scale real-time
  ride-sharing systems to address congestion issues faced by many
  cities. The goal is to serve all customers (service guarantees) with
  a small number of vehicles while minimizing waiting times under
  constraints on ride duration. This paper proposes an end-to-end
  approach that tightly integrates a state-of-the-art dispatching
  algorithm, a machine-learning model to predict zone-to-zone demand
  over time, and a model predictive control optimization to relocate
  idle vehicles. Experiments using historic taxi trips in New York
  City indicate that this integration decreases average waiting times
  by about 30\% over all test cases and reaches close to 55\% on the
  largest instances for high-demand zones.
\end{abstract}

%% file: introduction.tex
\section{Introduction}
\label{sec:intro}

Transportation Network Companies (TNCs) like Uber
and Lyft have fundamentally transformed mobility in many cities,
providing on-demand door-to-door transportation through mobile
applications.  They have also increased traffic in many cities: a
recent study by \citep{Erhardt2019} showed that, between 2010 and
2016, weekday vehicle hours of traffic delay have increased by 62\% in
San Francisco. In contrast, it was estimated that the delay increase
would have been 22\% without TNCs.  To address this issue, several
cities have begun limiting the number of TNC vehicles on the
road. Another way to tackle the underlying congestion and pollution
issues is to build mobility systems that utilize ride-sharing
\mbox{systematically}.  A study by \citep{Alonso-Mora462} showed that
systematic ride-sharing may significantly reduce the number of
vehicles needed to serve requests.  Their results indicate that 98\%
of the historic demand for taxi services in NYC could be served with a
much smaller taxi fleet, while maintaining short wait times.  This
paper continues this line of research and focuses on how to build a
real-time dispatching and routing architecture that serves the needs
of large-scale ride-sharing systems. It is envisioned that, in the
future, these ride-sharing systems will be deployed using autonomous
vehicles, guarantee service to all customers, and leverage advanced AI
systems. In the transition period, they can be supported by human
drivers provided that these drivers follow the instructions of the
ride-sharing system.

The value of stochastic information in real-time vehicle routing has
been demonstrated previously by \citep{scenariopvh}.  However,
incorporating stochastic information in large-scale ride-sharing
systems, where requests may arrive every tenth of a second during peak
times, is a challenge.  It is thus not surprising that
state-of-the-art approaches are purely myopic
\cite{Alonso-Mora462,riley2019}.  These systems batch requests and
optimize frequently to account for real-time demand.  Other approaches
to real-time dispatching (e.g., \cite{Holler2019}) use deep
reinforcement learning, but they ignore ride-sharing and do not
leverage advanced routing algorithms, focusing only on customer
assignment. To incorporate stochastic information, this paper
proposes a novel end-to-end framework (\name{}) for real-time
ride-sharing systems that tightly integrates state-of-the-art
optimization techniques, machine learning, and model predictive
control.

\begin{figure}[!t]
    \centering
    \includegraphics[width=\columnwidth]{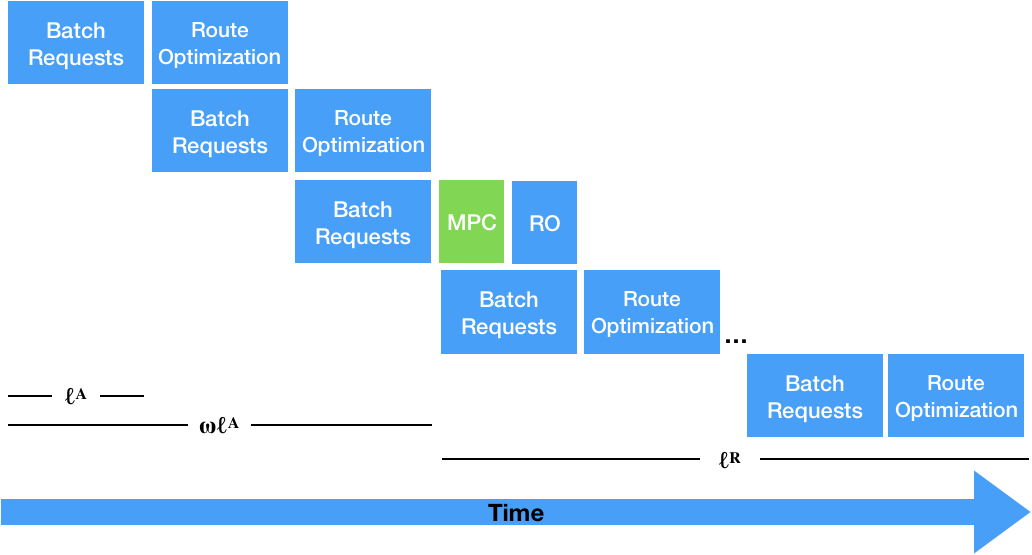}
    \caption{The \name{} Architecture for Real-Time Dispatching.}
    \label{fig:mpc_illustrated}
\end{figure}

The \name{} architecture is illustrated in
Figure~\ref{fig:mpc_illustrated}. Time is divided into epochs and,
during epoch $t$, \name{} optimizes the routing of the requests that
were batched in epoch $t-1$ as well as unserved requests from earlier
epochs. Moreover, at a lower frequency and prior to the routing
optimization, \name{} relocates idle vehicles using a Model Predictive
Control (MPC) step. The MPC step does not operate on individual
requests for scalability reasons. Rather it works with longer time
periods and at a coarser zone level (e.g., taxi zones in New York
City), and relies on a machine-learning model to predict the number of
requests between each pair of zones over time.  \emph{The main
  contribution of \name{} is to demonstrate that, in large-scale
  real-time ride-sharing systems, hybridizing state-of-the-art
  optimization algorithms for fine-grained routing decisions with
  model predictive control for idle vehicle relocation at a coarser
  space and time granularity provides significant operational
  benefits.}  Indeed, results on historic taxi trips from the New York
City Taxi and Limousine Commission \cite{nycdata}, indicate that this
tight integration decreases average waiting times by about 30\% over
all test cases and reaches close to a 32\% reduction in average
waiting times for high-demand zones in the most challenging instances.

This paper is organized as follows.
Section~\ref{sec:problem-statement} presents the problem.
Section~\ref{sec:related-work} presents related
work. Section~\ref{section:overall} gives an overview of \name{}.
Section~\ref{sec:rtdar} describes the dial-a-ride optimization.
Sections~\ref{sec:prediction} and~\ref{sec:idle-vehicle-relocation}
present the core conceptual contributions of the paper: the demand
forecasting model and the vehicle relocation
scheme. Section~\ref{sec:results} reports the experimental results and
Section~\ref{sec:conclusion} concludes the paper.

\section{Problem Statement}
\label{sec:problem-statement}

Operating a real-time ride-sharing system requires solving large-scale
dial-a-ride problems, where each request corresponds to a trip for a
number of riders from an origin to a destination that must take place
after a specified pickup time.  Constraints limit the time a rider can
spend inside a vehicle (ride duration constraints) and the number of
riders in a vehicle at any one time (vehicle capacity
constraints). The goal is to serve all requests and minimize the
average waiting time, while satisfying the ride duration and capacity
constraints.  Special attention is also devoted to ensure that no
request is left unserved indefinitely.  The systems studied in this
paper either use a fleet of autonomous vehicles or their own pool of
drivers who follow routing instructions exactly.  The system can thus
relocate the vehicles at will in order to anticipate demand.  It is
assumed that significant historical data is available and can be used
to forecast demand at the zone level.

%% file: related-work.tex
\section{Related Work} 
\label{sec:related-work}

A comprehensive review of popular dial-a-ride formulations which serve
as the foundation of \name{} can be found in \cite{Cordeau2007}.
\name{} uses a rolling horizon, alternating request batching and
optimization, as traditionally used in taxi pooling
\cite{stars,scalable-taxi}.  \citep{Alonso-Mora462} were first to
demonstrate the value of ride-sharing in NYC: they showed that 98\% of
the historic demand could be served with a smaller taxi fleet and
short wait times.  Their anytime algorithm uses cliques to generate
vehicle routes and hard-time windows to discard requests which cannot
be served efficiently.  A linear program is employed to move idle
vehicles towards discarded requests in order to better serve those
areas in the future. \citep{ZAC2018} improve over
\cite{Alonso-Mora462} by partitioning the region into zones and
assigning vehicles to zone paths. \citep{riley2019} is the first
algorithm designed to serve {\em all requests} with small waiting
times: they use column generation to serve all requests with smaller
number of vehicles and shorter waiting times. Their dial-a-ride
algorithm is used as the dispatching engine of \name{}. {\em To the
  author's knowledge, no algorithm other than the one of
  \citep{riley2019} provides guarantees to serve all requests: They
  can decide to ignore arbitrary requests.} Note that these three
algorithms are myopic: they do not exploit information about future
requests. \citep{Iglesias2017} proposed a model predictive control
approach for serving individual requests at the zone level, combining
a machine-learning model (based on deep learning) and a mixed-integer
program for request assignments and vehicle relocation. They did not
consider ride-sharing and their dispatching algorithm is performed at
a coarse granularity. This paper leverages and generalizes their model
predictive control approach. \citep{Ma2019} integrated dispatching
optimization and model predictive control for scheduling requests in a
multimodal transit systems: they do not batch requests, use a single
period for vehicle relocation, and assume Poisson arrivals for each
zone. The dispatching of each request uses local search and insertion
heuristics. The benefits of demand prediction and vehicle relocation
has been demonstrated by Bent and Van Hentenryck
\shortcite{Bent2007,scenariopvh} for various types of vehicle routing
problems (using online stochastic optimization) and by \citep{xyu2019}
for on-demand ride-pooling, using approximate dynamic
programming. \citep{Holler2019} used deep learning and bipartite
matching for dispatching and vehicle relocation: their approach does
not support ride sharing. \citep{Shah2020} enhance
\cite{Alonso-Mora462} with an approximation of the future reward
learned using a deep neural network. They provide improvements over
\cite{Alonso-Mora462} when the ride duration is twice as long as the
shortest path. However, the approach does not provide service
guarantees and does not minimize waiting times. It also rejects
requests even when vehicles are available, which can be problematic to
justify in practice. In contrast, this paper serves all requests with
an average waiting time of 2.5 minutes with 2,000 vehicles during peak
times and a more realistic ride-duration constraint (50\%
increase). The socio-economic benefits of ride-sharing systems is
explored by \citep{Bistaffa2019}. To the authors' knowledge, this
paper is the first integration of advanced optimization techniques,
machine learning, and model predictive control for the real-time
vehicle dispatching and relocation of large-scale ride-sharing
systems.

%% file: rtdar.tex
\section{Overall Organization}
\label{section:overall}

This section gives an overview of the \name{} architecture. As
depicted in Figure \ref{fig:mpc_illustrated}, \name{} divides time
into epochs of length $\ell^A$, i.e.,
$[0,\ell^A),[\ell^A,2\ell^A),[2\ell^A,3\ell^A),\ldots$. During epoch
      $\tau$, \name{} batches incoming requests and performs an
      optimization that assigns prior requests to vehicles and routes
      them. The requests considered in this optimization are those
      batched in epoch $\tau-1$, as well as unserved requests from
      earlier epochs. Periodically, \name{} performs a relocation
      optimization, which exploits a forecasting model to direct idle
      vehicles towards expected demand.

\paragraph{The Optimization Problem}

The optimization problem receives as inputs a set of requests, each of
which is characterized by its origin and destination, its earliest
pickup time, and its number of riders. The optimization has at its
disposal a number of vehicles. Each vehicle is characterized by its
departure location, its earliest departure time, its capacity, and its
set of riders. Each rider is characterized by her dropoff location and
the time she has already spent in the vehicle.

The starting location and departure time of a vehicle are given by the
current state of the mobility system. If a vehicle is idle in the
existing schedule, its starting location is its current position and
its departure time is the beginning of the epoch (i.e., $\tau
\ell^A$). If a vehicle is serving customers, its starting location is
the first location it visits after the start of the epoch and the
departure time is specified accordingly. The riders associated with a
vehicle are those who have already been picked up and need to be
dropped off. Hence, for every epoch, the optimization problem
considers all the requests whose riders have not been picked up yet,
while also scheduling the dropoffs of existing riders. Note that the
optimization problems associated with two successive epochs may
schedule a request differently as long as the request's riders have not
been picked up. This gives a lot of flexibility to the real-time
system at the cost of more complex optimization problems.

Given the computational complexity of the dial-a-ride problem that
must be solved in real time, the optimization may not be able to serve
all requests for some epochs. Hence, following \citep{riley2019},
\name{} associates a penalty with each request to ensure that the
request is served in reasonable time. The penalty is increased after
each epoch in which the request is not served. The optimization model
minimizes a weighted sum of the average waiting time and the penalties
associated with unserved requests.

\paragraph{Vehicle Relocation} Every $\omega$ epochs, \name{} performs
a relocation of vehicles at the zonal level (e.g., taxi zones, census
tracks, or traffic analysis zones).  The goal is to determine the
desired number of idle vehicles to move from zone $i$ to zone $j$ over
the next period $(\tau \ell^A, \tau \ell^A + \ell^R)$, where $\ell^R$
is the length of the relocation period and is significantly larger
than the epoch length $\ell^A$. As a result, the relocation
optimization operates at a much coarser granularity both in space and
time.

\emph{This combination of micro- and macro-decisions for routing and
  relocation is one of the salient features of \name{}} and is driven
by the reality of the large-scale real-time ride-sharing systems,
where the number of requests in each epoch makes it difficult
computationally to exploit forecasting information during the routing
decisions.

\paragraph{Forecasting the Demand} To inform the vehicle relocation, \name{}
is assumed to have at its disposal historical data for the number of
requests from zone $i$ to zone $j$ for every time period of length
$\ell^R$. This historical data is used to train a forecasting model of
the demand.

\section{The Dial-A-Ride Optimization}
\label{sec:rtdar}
\label{section-cg}

During each epoch, \name{} solves a generalized dial-a-ride
optimization specified in Section \ref{section:overall}. To perform
this task, it borrows the algorithm from \citep{riley2019}, which is
briefly reviewed in this section. The dial-a-ride optimization is
based on a column generation that operates at the route level. A
vehicle route specifies a sequence of pickups and dropoffs which
satisfies the ride duration constraints and the vehicle capacity. The
column generation interleaves the solving of (the linear relaxation
of) a Restricted Master Problem (RMP), which selects routes, and
pricing subproblems which generate new routes for each vehicle. The
process terminates when no new routes can improve the solution of the
RMP or the time limit for the column generation is met. The last stage
of the dial-a-ride optimization is a mixed-integer program that solves
the RMP exactly for the generated routes. The pricing subproblems are
complex due to their objective of minimizing waiting times. As a
result, traditional dynamic programming formulations are not effective
and \citep{riley2019} use an anytime exact algorithm that generates
routes of increasing lengths.

The RMP is depicted in Figure \ref{fig:master}. In the formulation,
$R$ denotes the set of routes, $R_v$ denotes the subset of possible
routes for vehicle $v$, $c_r$ represents the wait times incurred by
all customers served by route $r$, $p_{i}$ is the penalty of not
scheduling request $i$ for this epoch, and $a_{i}^{r}=1$ iff request
$i$ is served by route $r$. The RMP uses the following decision
variables: $y_r \in [0,1]$ is 1 iff route $r$ is selected and $z_{i}
\in [0,1]$ is 1 iff request $i$ is not served by any
of the selected routes.  The objective function minimizes the waiting
times of the served customers and the penalties for the unserved
customers.  Constraints~\eqref{model:master_constr:allserved} ensure
that $z_{i}$ is set to 1 if request $i$ is not served by any of the
selected routes and constraints~\eqref{model:master_constr:oneroute}
ensure that only one route is selected per vehicle.

\begin{figure}[!t]
\begin{subequations} \label{model:master}
\begin{align}
\min 	&  \quad \sum_{r \in R}  c_r y_r + \sum_{i \in P} p_{i} z_{i} \label{model:master_obj}\\
\mbox{subject to} &  \left ( \sum_{r \in R} y_r a_{i}^{r} \right ) + z_{i} = 1 & \forall i \in P \label{model:master_constr:allserved} \\ 
	&  \sum_{r \in R_v} y_r = 1 & \forall v \in V \label{model:master_constr:oneroute} \\
	&  z_{i} \in \mathbb{N} & \forall i \in P \label{model:master_constr:domainz}\\
	&  y_r \in \{0,1\} & \forall r \in R \label{model:master_constr:domainy}
\end{align}
\end{subequations}
\caption{The Resricted Master Problem Formulation.}
\label{fig:master}
\end{figure}

Since the dial-a-ride optimization may not schedule all the requests,
it is important to update the penalty of unserved requests to ensure
that they will not be delayed too long. For the penalty for an
unserved request $c$ in epoch $\tau$, \citep{riley2019} use $p_c =
\rho 2^{(\tau \ell^A - e_c) / (10\ell^A)} \label{equation:obj} $,
where $e_c$ is the earliest possible pickup time for request $c$.  The
$\rho$ parameter was tuned to incentivize the algorithm to schedule
each incoming request in its first available epoch.

%% file: prediction.tex
\section{Demand Forecasting}
\label{sec:prediction}


Forecasting the demand from zone $i$ to zone $j$ over time may be
challenging in some settings, since this demand may be sparse for some
zones and historical data may be limited. To address this difficulty,
\name{} proceeds in two steps: it first predicts the number of
requests in a zone $z$ in time period $t$ and then approximates the
zone-to-zone demand.

\paragraph{Preprocessing}

Let $d_{zt}$ be the demand for zone $z$ during period $t$. In the case
study, the time series $\{d_{zt}\}_{t}$ is strongly non-stationary
(the mean and variance vary over time). As a result, the forecasting
model first stationarizes the time series by differencing it over a
week period. More precisely, the forecasting model defines
$\delta_{zt} = d_{zt} - d_{z(t-n_r\times7)}$, where $n_r$ is the
number of periods in a day, and predicts the differenced demand
$\delta_{zt}$ instead of $d_{zt}$.

\paragraph{Vector Autoregression}

To forecast the time series $\{\delta_{zt}\}_{t}$, \name{} uses
Vector Autoregression (VAR), a multivariate generalization of
Autoregression (AR). In VAR, the expected value of a multivariate time
series at a particular period is assumed to be a linear function of
the value of the time series at previous time steps.

The prediction for the differenced demand in zone $z$ in period $t$
uses not only $z$'s demand in prior periods but also the differenced
demands of $z$'s adjacent zones. Let $N(z)$ denote the zones adjacent to $z$
and $d=|N(z)|+1$. Let vector $\Delta_{zt}\in \mathbb{R}^d$ denote the
weekly-differenced demands of $z$ and its adjacent zones in period
$t$. Each element in $\Delta_{zt}$ is an element in
$\{\delta_{zt}\}_{z\in N(z)\cup \{z\}}$. $\delta_{zt}$
can then be modeled as
\begin{equation}\label{var0}
    \delta_{zt} = \phi_{zt-1} \Delta_{zt-1} + \dots + \phi_{zt-k} \Delta_{zt-k} + \eta
\end{equation}
where $\phi_{zt}$ is a row vector in $\mathbb{R}^d$ and $\eta$ is a
white noise with zero mean. The coefficients $\phi_{zt}$ are estimated
using least square regression and the order $k$ is selected based on
the \textit{Akaike information criterion (AIC)}.

Once the parameters have been estimated, the prediction for the
differenced demand of $z$ in period $t$ is given by
\begin{equation*}\label{var}
    \bar{\delta}_{zt} = \bar{\phi}_{zt-1} \Delta_{zt-1} + \dots + \bar{\phi}_{zt-k} \Delta_{zt-k}
\end{equation*}
\noindent
The demand prediction for zone $z$ at time $t$ is then given by $
\bar{\lambda}_{zt} = d_{z(t-n_r\times7)} + \bar{\delta}_{zt}.  $

\paragraph{Destination Assignment}

Given the number of requests for zone $z$ in period $t$, the trip
destinations for these requests are assigned using historical
distributions. \name{} uses an historical distribution for each hour
during the weekdays and the weekend days. For example, when predicting
the demand in zone $z$ during the 7--8am period on a Wednesday, if
$70$ percent of the trips originating from $z$ during this period on
weekdays have their destination in zone $b$ in historical data, then
the number of trips going from $z$ to $b$ will be
$0.7\bar{\lambda}_{zt}$ rounded to the nearest non-negative
integer. This returns the final demand prediction
$\bar{\lambda}_{ijt}$ for the requests from zone $i$ to zone $j$ at
$t$.

%% file: vehicle-relocation.tex
\section{Idle Vehicle Relocation}
\label{sec:idle-vehicle-relocation}

The idle vehicle relocation process is run every $\omega$ epochs and
considers periods of length $\ell^R$, i.e.,
$[0,\ell^R),[\ell^R,2\ell^R),\ldots$ It has at its disposal the zone
    to zone demand forecasts for each period. It proceeds in two
    steps: (1) It first uses a Model Predictive Control (MPC) approach
    to find the desired number of vehicles in each zone; (2) It then
    selects the vehicle relocation assignments to minimize the
    relocation cost.

\paragraph{Zone Rebalancing}
\label{sec:zone-balancing}

The MPC approach in this section is borrowed from \citep{Iglesias2017}
and generalized to real-time ride-sharing systems, where multiple
riders can share a vehicle. Its goal is to determine the number of
idle vehicles to move from zone $i$ to zone $j$ during the next period
$(\tau \ell^A, \tau \ell^A + \ell^R)$ in order to minimize the average
waiting time in the dial-a-ride optimizations. To achieve this goal,
the MPC approach solves an assignment optimization at the zone level
over multiple time periods. Hence, in contrast to the dial-a-ride
optimization, it works at a coarser granularity and looks ahead in the
future using the demand forecasting module.

The MPC approach uses a MIP model (MPC-MIP) over $T$ periods, each of
length $\ell^R$. Let $\mathcal{T} = \{0, \dots, T-1\}$. For each
period $t \in \mathcal{T}$, MPC-MIP takes as input
$\bar{\lambda}_{ijt}$, the forecasted demand originating in zone $i$
with a destination in zone $j$ at time $t$, as well as a variety of
information about the current state of the system. In particular,
$w_{ij}$ is the current sharing ratio for requests from zone $i$ to
zone $j$ in the system (e.g., $w_{ij} = 1$ means that riders are alone
in the vehicle, while $w_{ij} = 1.5$ means an average of $1.5$
passengers per vehicle); $A_{it}$ is the set of vehicles that will
become idle in zone $i$ during period $t$ (estimated from routes of
current vehicles) and $tt_{ij}$ is the average number of periods it
takes to move from a stop in zone $i$ to a stop in zone $j$.

The MPC-MIP decision variables are: the number $x^r_{ijt}$ of empty
vehicles to move from zone $i$ to zone $j$ starting at period $t$; the
number $x^p_{ijt}$ of vehicles with passengers moving from $i$ to $j$
starting at period $t$; and the number $u_{ijt}$ of requests
originating in zone $i$ and ending in zone $j$ which are not served at
period $t$. The objective~\eqref{mpc:obj} minimizes the number of unserved
requests while also enforcing a small penalty on moving vehicles.
This penalty ensures that vehicles prefer to stay in their current
zone if that current zone is expected to need them in the
future.  Constraints~\eqref{mpc:constr:demand_balance} and
\eqref{mpc:constr:demand_balance} are the flow balance constraints for
requests.  Constraints~\eqref{mpc:constr:vehicles_balance} are the
flow conservation constraints for vehicles.

\begin{figure}[!t]
\begin{subequations} \label{mpc}
\begin{flalign}
& \text{min} \sum_{t=0}^{T-1}\sum_{i,j \in Z} (T - t)u_{ijt}  + \sum_{t=0}^{T-1} \sum_{i,j \in Z} tt_{ij} x^r_{ijt} \label{mpc:obj} \\
\intertext{subject to} 
& \; x^p_{ij0}  + u_{ij0} = \ceil{\bar{\lambda}_{ij0} / w_{ij}} \quad (\forall i,j) \label{mpc:constr:demand_balance_0} \\	
& \; x^p_{ijt}  + u_{ijt} = \ceil{\bar{\lambda}_{ijt} / w_{ij}} + u_{ijt-1} \quad (\forall i,j,t) \label{mpc:constr:demand_balance} \\ 
& \textstyle \; \sum_{j} x^p_{ijt} + x^r_{ijt} = |A_{it-1}| + \sum_{j} x^p_{jit - tt_{ji}} + x^r_{jit - tt_{ji}} \; (\forall i,t)  \label{mpc:constr:vehicles_balance} \\
& \; x^p_{ijt} \in \mathbb{Z}, x^r_{ijt} \in \mathbb{Z}, u_{ijt} \in \mathbb{Z}  \label{mpc:constr:domain}
\end{flalign}
\end{subequations}
\caption{The MPC-MIP Model for Zone Balancing.}
\label{fig:mpc}
\end{figure}


\paragraph{Vehicle Relocation}
\label{sec:vehicle-relocation}

MPC-MIP returns the number $\bar{x}^r_{ijt}$ of vehicles that should
move from zone $i$ to zone $j$ in period $t$. Only the relocations in
period $\bar{x}^r_{ij0}$ are relevant for \name{} at this point in
time.  However, \name{} must now identify $\bar{x}^r_{ij0}$ specific
vehicles to relocate. This is performed by another MIP model (VR-MIP),
which receives the following inputs: $\bar{x}^r_{ij0}$, the sets
$A_{i0}$ of idle vehicles in zone $i$, the time $c_{vj}$ to move
vehicle $v$ to its closest stop in zone $j$.
VR-MIP decides whether a vehicle is relocated to a zone: Variable
$y_{vj}$ is 1 if vehicle $v$ is chosen to relocate to zone $j$.  The
VR-MIP objective~\eqref{mpc-2:objective} minimizes the sum of the
traveling times, Constraints~\eqref{mpc-2:constr:move_correct_vehicles}
ensure the correct number of relocations from zone $i$ to zone $j$,
and Constraints~\eqref{mpc-2:constr:only_move_to_one_zone} ensure that
a vehicle does not relocate to more than one zone.

\begin{figure}[!t]
\begin{subequations} \label{mpc-2}
\begin{align}
\text{min} \;\;\;\;\;	& \sum_{v \in A_{i0}} \sum_{j} c_{vj} y_{vj} \label{mpc-2:objective} \\
\mbox{subject to}
& \sum_{v \in A_{i0}} y_{vj} = \bar{x}^r_{ij0} & \forall j \in Z \label{mpc-2:constr:move_correct_vehicles} \\
& \sum_{j} y_{vj} \leq 1 & \forall v \in A_{i0} \label{mpc-2:constr:only_move_to_one_zone} \\
& y_{vj} \in \{0, 1\} & \forall v \in A_{i0}, j \in Z
\end{align}
\end{subequations}
\caption{The Vehicle Relocation Optimization (VR-MIP).}
\label{fig:mpc-2}
\end{figure}

%% file: results.tex
\section{Experimental Results} \label{sec:results}

\paragraph{Case Study}

\name{} is evaluated on the yellow trip data obtained via the New York
City Taxi and Limousine Commission (NYCTLC) \cite{nycdata}.  The
NYCTLC dataset provides the number of passengers for each trip and the
start time of each trip, which is used as both the request time and
the lower bound on the earliest pickup time. This section reports
results on 24 instances, two hours each day for two days per month
from July 2015 through June 2016. Rush hours (7--9am) were selected to
obtain instances that are computationally challenging. The instances
have an average of \num[group-separator={,}]{48100.5} customers and
range from \num[group-separator={,}]{19276} to
\num[group-separator={,}]{59820} customers. Individual requests with
more customers than the vehicle capacity are split into several
trips. Following \citep{riley2019}, in order to ease ride sharing,
avoid curb management issues, and reduce the number of stops,
Manhattan is represented by a grid with cells of 200 squared
meters. Each such cell represents a pickup/dropoff location. The
travel time matrix for the network of locations was then precomputed
by querying OpenStreetMap \cite{OpenStreetMap}. For every trip, the
locations of the origin and destination were obtained by selecting the
closest locations to their pickup and dropoff points in the NYCTLC
dataset.

\paragraph{Runtime Configurations}

\name{} is compared to its myopic version \oname{} which has no idle
vehicle relocation and is essentially the approach proposed by
\citep{riley2019}. It is also compared to \oracle{}, a version of
\name{} using perfect information on future requests instead of the
machine-learning predictions. Unless otherwise specified, all
experiments were performed with the following default parameters for
\name{} (and \oname{} when relevant): \num[group-separator={,}]{2000}
vehicles of capacity 4, a maximum deviation from the shortest path
determined by $\max \{ \alpha t_c, \beta + t_c \}$, where $t_c$ is the
shortest possible path from customer $c$'s origin to their
destination, $w_{ij} = 1.2$, $\rho = 420$, $\alpha = 1.5$, $\beta =
240$ seconds, $\ell^A = 30$ seconds, $\ell^R = 300$ seconds, and
$\omega = 10$ epochs.  Empty vehicles are initially distributed evenly
over the locations.  The demand was forecasted at the hour level and
scaled down uniformly due to the sparse demand in some of the
zones. All MIP models are solved using Gurobi~8.1 \cite{gurobi}.

\paragraph{Reduction in Waiting Times}

\begin{figure}[!t]
  \centering
  \includegraphics[width=0.75\columnwidth]{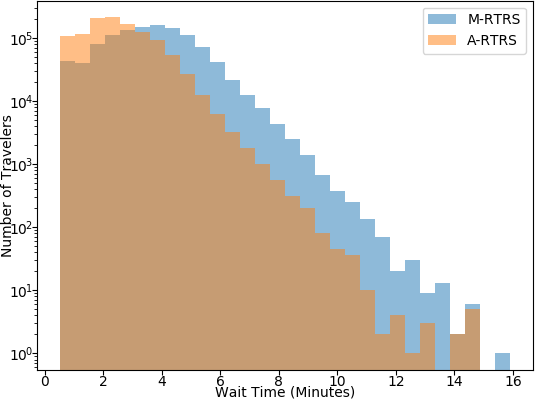}
  \caption{Histograms of Waiting Times for \name{} and \oname{} (Logarithmic Y-Scale).}
  \label{fig:hist}
\end{figure}

\begin{figure}[!t]
  \centering
  \includegraphics[width=0.75\columnwidth]{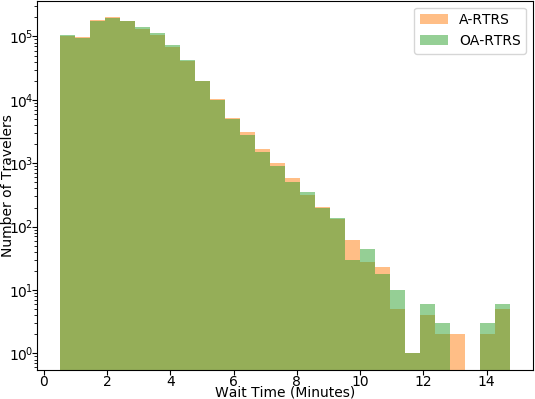}
  \caption{Histograms of Waiting Times for \name{} and \oracle (Logarithmic Y-Scale).}
  \label{fig:oracle_hist}
\end{figure}

Figures~\ref{fig:hist} and \ref{fig:oracle_hist} report the distributions of the waiting times
incurred by all customers across all instances with a logarithmic
y-axis. The results demonstrate that \name{} reduces waiting times
across the vast majority of trips. It reduces the average waiting times
from 3.64 to 2.51 minutes (a 30\% improvement) while also decreasing
the standard deviation from 1.48 to 1.16.

\begin{figure}[!t]
  \centering
  \includegraphics[width=0.75\columnwidth]{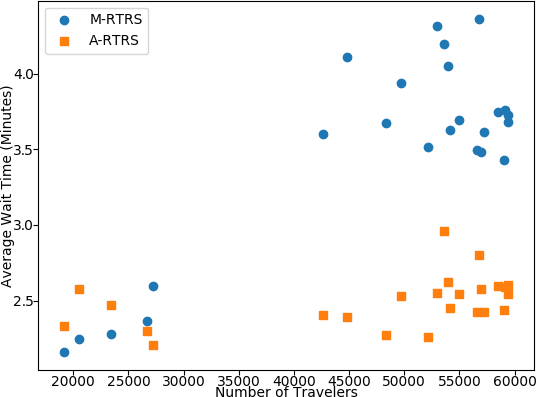}
  \caption{Waiting Times for \name{} and \oname{}.}
  \label{fig:cpaior_compare}
\end{figure}

\begin{figure}[!t]
  \centering
  \includegraphics[width=0.75\columnwidth]{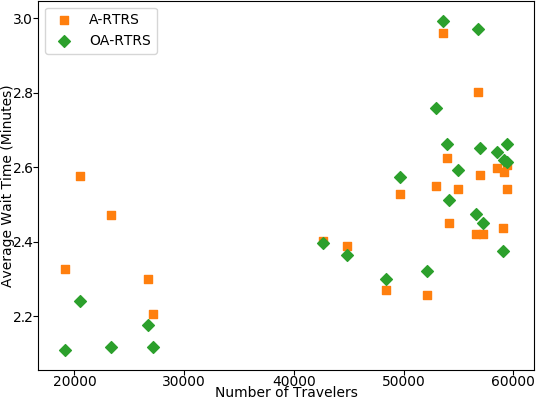}
 \caption{Waiting Times for \name{} and \oracle{}.}
  \label{fig:oracle_compare}
\end{figure}

Figures~\ref{fig:cpaior_compare}, \ref{fig:oracle_compare}, and
\ref{table:avg_wait} go into more details and report results on each
instance and by instance sizes. The results show that the benefits of
relocation are more significant for instances with a large number of
requests. \name{} strongly dominates \oname{} for instances with large
demands, but is not as effective on instances with relatively low
demand. This last result is due to the accuracy of the
machine-learning algorithm as highlighted in Figure
\ref{fig:oracle_compare}. With perfect information, \oracle{} improves
over \name{} over low-demand instances, but \name{} and \oracle{}
behave similarly on high-demand instances.

\begin{figure}[!t]
\begin{center}
 \begin{tabular}{||c | c c c||} 
 \hline
 Number of Requests & \oname{} & \name{} & OA-RTRS \\ [0.5ex] 
 \hline\hline
 $< \num[group-separator={,}]{40000}$ & 2.33 & 2.37 & 2.15 \\ 
 \hline
 $\num[group-separator={,}]{40000}-\num[group-separator={,}]{50000}$ & 3.83 & 2.40 & 2.41 \\
 \hline
 $\num[group-separator={,}]{50000} < $ & 3.78 & 2.56 & 2.56 \\
 \hline
\end{tabular}
\end{center}
\caption{Average Waiting Times by Instance Sizes.}
\label{table:avg_wait}
\end{figure}

\begin{figure}[!t]
  \centering
 \includegraphics[width=\columnwidth]{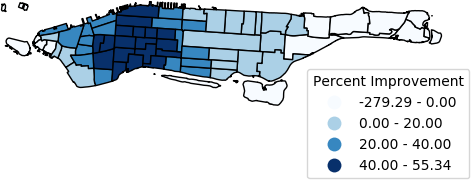}
  \caption{Percentage Improvement in Waiting Times by Origin.}
  \label{fig:zone_wait_map}
\end{figure}

\begin{figure}[!t]
  \centering
  \includegraphics[width=\columnwidth]{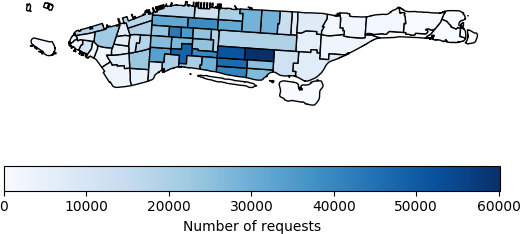}
  \caption{Total Number of Requests By Zone (over all instances).}
  \label{fig:requests_per_zone}
\end{figure}

\begin{figure}[!t]
  \centering
  \includegraphics[width=\columnwidth]{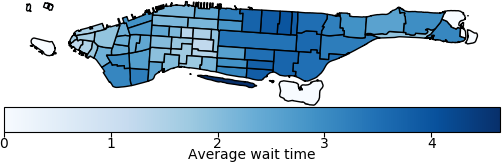}
  \caption{Average Waiting Times by Zone.}
  \label{fig:wait_per_zone}
\end{figure}

\paragraph{Zonal Information} Figures~\ref{fig:zone_wait_map}~and~\ref{fig:requests_per_zone}
describe where the reductions in waiting times occur. Together, they
show that relocation benefits most the zones with significant demand,
where the improvements can reach almost 55\%. Figure
\ref{fig:wait_per_zone} is reassuring from a fairness standpoint:
Zones with low demand keep low average waiting times under relocation.

\begin{figure}[!t]
  \centering
  \includegraphics[width=0.75\columnwidth]{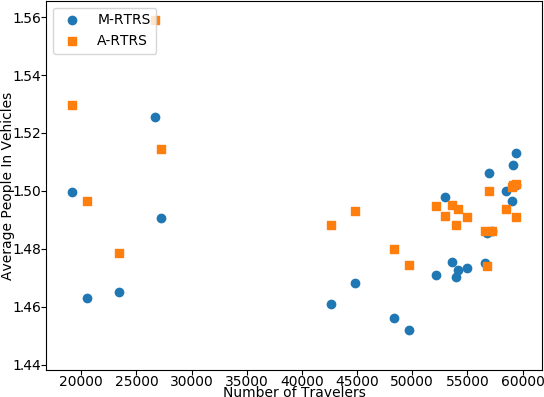}
  \caption{Average Number of Passengers per Vehicle.}
  \label{fig:cpaior_vehicle_util_no_deadhead}
\end{figure}

\paragraph{Passenger Information} Figures~\ref{fig:cpaior_vehicle_util_no_deadhead} and
\ref{fig:oracle_vehicle_util_no_deadhead} show that relocation
slightly decreases ride sharing which is beneficial for
passengers. There is less of a need to ride share to satisfy the
demand and minimize waiting times.

\begin{figure}[!t]
  \centering
  \caption{Average Number of Passengers per Vehicle.}
  \label{fig:oracle_vehicle_util_no_deadhead}
\end{figure}

\begin{figure}[!t]
    \centering
  \includegraphics[width=0.75\columnwidth]{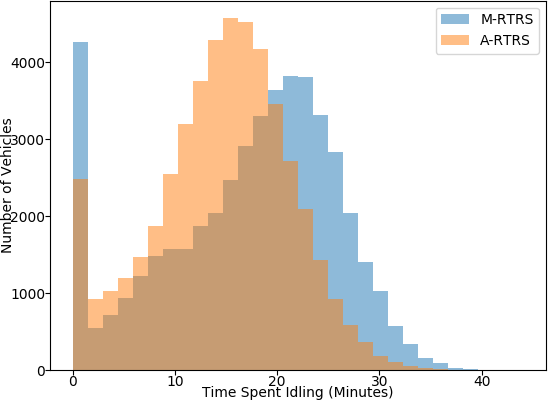}
    \caption{Idling Times per Vehicle.}
    \label{fig:cpaior_idling_hist}
\end{figure}

\paragraph{Vehicle Information} Figures~\ref{fig:cpaior_idling_hist}
depicts the histograms of idle times per vehicles. It shows that \oname{}
has more vehicles with no idle time and more vehicles with high idle times.
\name{} is more balanced. 

\paragraph{Relocation} 

\begin{figure}[!t]
    \centering
    \includegraphics[width=0.75\columnwidth]{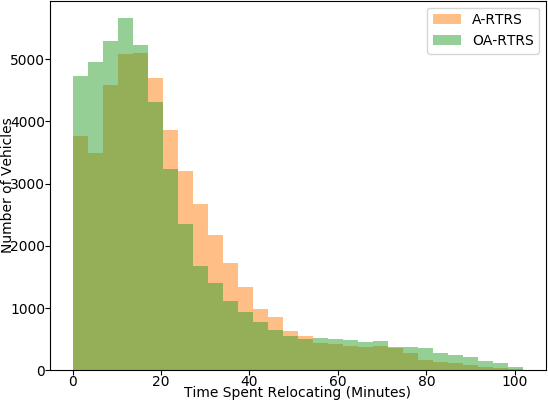}
    \caption{Relocation Time per Vehicle.}
    \label{fig:relocating_hist}
\end{figure}

\begin{figure}[!t]
   \centering
    \includegraphics[width=0.75\columnwidth]{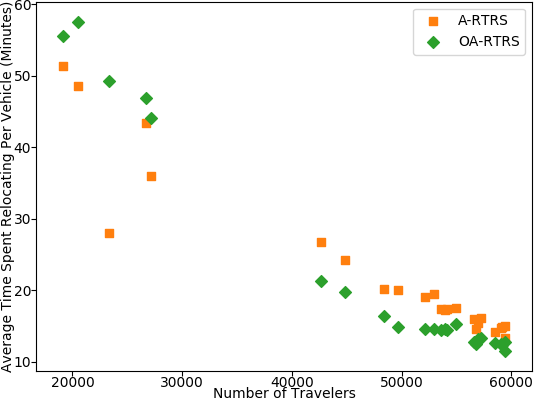}
    \caption{Average Times Vehicles Spend Relocating.}
    \label{fig:relocating_scatter}
\end{figure}

Figure ~\ref{fig:relocating_hist} shows that most vehicles spend less
than five minutes in relocation and the vast majority of vehicles
spend less than 17\% of their operating hours relocating. Figure
\ref{fig:relocating_scatter} shows that \oracle{} relocates more on
the instances with lower demands, explaining why it improves over
\name{} on these instances. Having a better demand predictor is thus
an important research direction.

%% file: conclusion.tex
\section{Conclusion}
\label{sec:conclusion}

This paper proposed \name{}, an end-to-end framework for real-time
optimization of large-scale ride-sharing systems. \name{} combines
demand forecasting, state-of-the-art optimization, and model
predictive control to dispatch, route, and relocate vehicles in
real-time, minimizing average waiting times. The mobility system
provides service guarantees (i.e, it serves all requests), enforces a
ride-duration constraint (i.e., no passenger travels more than 50\%
over their shortest path), minimizes waiting times, while achieving
reasonable waiting times through penalties increasing over time.
Experiments using historic taxi trips in New York City indicate that
this integration decreases average waiting times by about 30\% over
all test cases and reaches close to 55\% on the largest instances for
high-demand zones compared to a base line without relocation.  On the
NYC case study, \name{} serves all requests in reasonable time and
with an average waiting of 2.51 minutes with a standard deviation of
1.16, using a fleet of 2,000 vehicles of capacity 4. The results also
demonstrate that, while zones with large demand see the most benefits,
zones with low demand maintain low waiting times and that the vast
majority of vehicles spend less than 17\% of their operating time
relocating. In summary, \name{} demonstrates that, in large-scale
real-time ride-sharing systems, hybridizing state-of-the-art
optimization algorithms for fine-grained routing decisions with model
predictive control for idle vehicle relocation at a coarser space and
time granularity provides significant operational benefits. Future
research will be devoted in improving the machine-learning algorithm,
since more accurate predictions will enable a better performance on
instances with relatively fewer requests.